\documentclass{conm-p-l}
\usepackage {amsmath, amscd, amsthm}

\newtheorem{thm}{Theorem}[section]

\newtheorem{lem}[thm]{Lemma}

\theoremstyle{definition}
\newtheorem{defn}[thm]{Definition}

\theoremstyle{remark}

\numberwithin{thm}{section}
\newtheorem{lem/defn}[thm]{Lemma/Definition}

\numberwithin{equation}{section}

\newcommand{\ten}{\otimes}
\newcommand{\del}{\bigtriangleup}
\newcommand{\second}{\prime \prime }

\DeclareMathOperator{\End}{End}
\DeclareMathOperator{\Sym}{Sym}

\begin{document}
\title[Symmetric polynomials and $H_D$-quantum vertex algebras ]{Symmetric
    Polynomials and $H_D$-Quantum Vertex Algebras}
    \author{I. I. Anguelova}
    \address[Iana Anguelova] {Department of Mathematics,
University of Illinois, Urbana, IL 61801, USA }
    \email{anguelov@math.uiuc.edu}
   

       \begin{abstract}  In this paper we  use a  bicharacter
    construction  to define an $H_D$-quantum vertex algebra structure  corresponding to the quantum vertex operators describing classes of  symmetric polynomials.  
\end{abstract}

\maketitle

\vspace{0.5cm}

\section{Introduction}
\label{sec:intro}
Vertex operators were introduced in the earliest days of string theory
and axioms for vertex algebras were developed to incorporate these
examples (see for instance \cite{FLM}).
Similarly, the definition of quantum vertex algebra should be such that it       accommodates  the existing examples of quantum vertex operators and their properties (see for instance \cite{FJ}, \cite{FR2}, \cite{BFJ}, \cite{JK} and many others).

In this  series of papers we study the quantum vertex algebra structure
corresponding to classes of symmetric polynomials (e.g., Hall-Littlewood
or  Macdonald polynomials). The vertex operators describing these
polynomials were considered by N. Jing in a series of papers
(\cite{J1}, \cite{J2}, \cite{J3}). As shown for instance by these examples a  major
difference between classical and quantum vertex algebras lies in the
fact that two  vertex operators (fields) $Y(z)$ and $Y(w)$ are no longer
`almost' commuting (i.e., commuting, except on the diagonal $z=w$), but instead there
is a braiding map connecting the products $Y(z)Y(w)$ and $Y(w)Y(z)$ (see section 3).  The goal
in this series of papers is to incorporate these vertex operators in
certain braided vertex algebra structures.

There are several proposals for the axioms of quantum vertex algebras
(or deformed  chiral algebras). We will consider throughout the papers
three of these definitions, namely by Borcherds in \cite{Borc}, Frenkel-Reshetikhin in \cite{FR} and Etingof-Kazhdan in \cite{EK}). We will refer to them  as $(A,H,S)$ quantum vertex algebras, deformed chiral algebras (of FR type) or quantum vertex algebras of EK type.
One of the goals of this series of papers is to show by examples that
these axioms are surprisingly not equivalent. Also we show
that the axioms of EK or FR are not sufficient to describe the vertex
operators of the symmetric polynomials.  

Our main goal  in this paper will be to define an $H_D$-type of
quantum vertex
algebra. (Here
$H_D=\mathbf{C}[D]$ is the Hopf algebra of infinitesimal translations,
a fundamental ingredient in the classical vertex algebras.) We
construct examples of such $H_D$-quantum vertex algebras incorporating
the quantum vertex operators appearing in the theory of
Hall-Littlewood and Macdonald symmetric polynomials. 

We use extensively the bicharacter construction from \cite{Borc} as a
help towards defining the quantum vertex algebra structures. We will
state  a technical theorem (see Theorem \ref{thm:main}) which will allow us to use this
construction on a large class of vertex operators corresponding to the
general orthogonal polynomials as defined by Macdonald.

\section{Symmetric polynomials and vertex operators}
\label{sec:sympol}
In this section we recall the Macdonald definition of general symmetric polynomials (\cite{Macd}). Also we recall the vertex operators associated to these general symmetric polynomials (as considered for instance by Jing, \cite{J1}, \cite{J3}) 

We work over a field $k$ of characteristic zero containing the rationals, $F$ is a field extension of $k$ (for instance $k(t)$ or $k(q,t)$, where $q,t$ are parameters).

Denote by  $\Lambda $ the ring of symmetric functions over $\mathbf{Z}$ in countably many independent variables $x_i,i \ge 0 $.
Let $\Lambda _{F}=\Lambda \ten_{\mathbf{Z}} F$.
As usual denote by $p_{n}$ the power symmetric functions.
If $\lambda$ is a partition, $\lambda =(\lambda _1, \lambda _2,
..,\lambda _k, ...),  \ \lambda _1 \ge \lambda _2 \ge  ...\ge \lambda
_k \ge ...$, we call a family $(a_{\lambda})$  of elements in a ring
indexed by partitions multiplicative if $a_{\lambda}=\prod a_{\lambda
  _i}$. 
The family $(p_{\lambda})$ is a multiplicative basis for the symmetric
functions. We will also use the basis $(m_{\lambda})$ of monomial
symmetric functions (it is not multiplicative).
Denote  $z_{\lambda}=\prod _{i\ge 0}i^{m_i}.m_i!$, where $m_i=m_i(\lambda)$ is the number of parts of $\lambda$ equal to $i$.

Let $(v_{\lambda})$ be a multiplicative family in $F^{\times }$. We want to
define for each such family a set of symmetric polynomials
$\{P_{\lambda}\}$ indexed by partitions. First, we define a scalar
product $\langle \ ,\ \rangle _{v}$ on $\Lambda _{F}$ (depending on the
multiplicative family $v_{\lambda }$) by 
\begin{equation*}
\label{eq:scalarpr}
\langle p_{\lambda}, p_{\mu} \rangle _{v} =\delta _{\lambda \mu} z_{\lambda}v_{\lambda},
\end{equation*}
for any partitions $\lambda ,\mu$.
The set $\{P_{\lambda}\}$ associated to $v_{\lambda }$ should satisfy  the following two (over-determining) conditions:
\begin{align*}
&\langle P_{\lambda}, P_{\mu} \rangle _{v}=0 \ \ \ \text{for} \ \ \lambda \neq \mu, \\
&P_{\lambda}=m_{\lambda}+\sum _{\mu < \lambda}u_{\lambda \mu}m_{\mu}, \ \ u_{\lambda \mu}\in F.
\end{align*}
Here $\mu < \lambda$ is with respect to the usual partial order on partitions.
In the cases where $\{P_{\lambda}\}$ exist denote also by
$Q_{\lambda}$ the dual of $P_{\lambda}$, i.e., $\langle
P_{\lambda},Q_{\mu}\rangle _{v}=\delta _{\lambda , \mu }$.

For general  $(v_{\lambda})$ the corresponding  $\{P_{\lambda}\}$ might
not exist, but the following important examples  are known:
\begin{itemize}
\item
If $v_n=1$ \  the family  $\{P_{\lambda}\}$ is in fact $\{S_{\lambda}\}$, the Schur symmetric functions; Schur functions are self-dual.
\item
If $v_n=\frac{1}{1-t^n}$ the   $\{P_{\lambda}\}$\  are  the Hall-Littlewood functions.
\item
If  $v_n=\frac{1-q^n}{1-t^n}$  \ the family $\{P_{\lambda}\}$ consists of the Macdonald symmetric functions.
\item 
Jack symmetric functions and the zonal polynomials can also be obtained by picking appropriate $v_n$, but we will  not consider them in this paper.
\end{itemize}

We can view $p_n$ as an operator acting on $\Lambda _F$ by multiplication. Define also the operators $p^{\perp}_{n}$ by requiring
\begin{displaymath}
\langle p^{\perp}_{n}f,g\rangle _{v}= \langle f,p_n g\rangle _{v},
\end{displaymath}
for any $f,g \in \Lambda _F$.

Given a multiplicative family $(v_{\lambda })$, let $\mathcal{H}=\mathcal{H_v}$ be the algebra generated by the operators 
$h_n=-v_{-n}^{-1}p^{\perp}_{-n}, \ h_{-n}=v_n ^{-1}p_n$ for $n\in \mathbf{N}$.
\begin{lem} 
\label{lem:heis}
For each multiplicative family $(v_{\lambda })$ the operators $\{ h_n | n\in \mathbf{Z}\}$ generate a representation of a deformed Heisenberg algebra on $\Lambda _F$, i.e.,
\begin{equation}
[h_m,h_n]=mv^{-1}_{|m|}\delta _{m+n,0}1.
\end{equation}
\end{lem}


Define the following vertex operator
\begin{equation}
\Phi(z)=\exp (\sum _{n\ge 1}\frac{h_{-n}}{n} z^n)\exp (-\sum _{n\ge 1}\frac{h_{n}}{n} z^{-n})=\sum \Phi _n z^n \ ,
\end{equation}
where $\Phi _n\in \End(\Lambda _F)$ are the Fourier coefficients of the field $\Phi (z)$.
(Note that in fact we have a different vertex operator for each multiplicative family $(v_{\lambda })$). 
The usefulness of these vertex operators is due to the following facts.
Let   $\lambda =(\lambda _1, \lambda _2, ..,\lambda _k)$ be a partition.
\begin{itemize}
\item
In the case $v_n=1$, 
\begin{displaymath}
\Phi _{\lambda _1}\Phi _{\lambda _2}...\Phi _{\lambda _k}1=S_{\lambda}.
\end{displaymath}
Here $S_{\lambda}$ is the Schur function corresponding to the partition $\lambda$ (recall that Schur functions are self-dual).
\item
In the case $v_n= \frac{1}{1-t^n}$, 
\begin{displaymath}
\Phi _{\lambda _1}\Phi _{\lambda _2}...\Phi _{\lambda _k}1=Q_{\lambda}.
\end{displaymath}
Here $Q_{\lambda}$ is the dual to the Hall-Littlewood polynomial corresponding to the partition $\lambda$. (This fact was proved by Jing in \cite{J1}).
\item
In the case of $v_n= \frac{1-q^n}{1-t^n}$ the symmetric function $\Phi _{\lambda _1}\Phi _{\lambda _2}...\Phi _{\lambda _k}1$ is connected to the Macdonald symmetric function, but not in general equal to the dual Macdonald   function as in the previous two cases. For a single row partition they are indeed equal, for two-row partitions Jing proves that the Macdonald polynomials are related via a $q$-hypergeometric function (see \cite{J5}).
\end{itemize}
Define  the space $V=\Lambda _F \ten \mathbf{C}[\mathbf{Z}\alpha]$, where $\mathbf{C}[\mathbf{Z}\alpha]$ is the group algebra of the rank-one lattice generated by $e^{\alpha}$, $m\in \mathbf{Z}$ (such that $e^{m\alpha}e^{n\alpha}=e^{(m+n)\alpha}$).
Define also the vertex operator on $V$
\begin{equation}
\label{eq:psi}
\Psi(z)=\exp (\sum _{n\ge 1}\frac{h_{-n}}{n} z^n)\exp (-\sum _{n\ge 1}\frac{h_{n}}{n} z^{-n})e^{\alpha}z^{\partial _{\alpha}},
\end{equation}
which is more convenient to work with.

In the case of Schur polynomials the vertex operator $\Psi (z)$ belongs to the  classical vertex algebra of the odd rank-one lattice (the Heisenberg algebra in this case is the undeformed infinite-dimensional Heisenberg algebra). 
In the deformed case the vertex algebra can no longer be classical,
i.e., the vertex operators will no longer ``almost'' commute, but
instead will form a ``braided singular ring'', i.e., they will belong to
a braided vertex algebra structure.  The goal of this paper will be to
describe the  quantum vertex algebra structure to which this vertex operator $\Psi (z)$ belong. 

The simplest case of a  braided vertex algebra structure---the  $H_D$-type---is described below.

\section{$H_D$-quantum vertex algebras }
\label{sec:EK type }

Let $V$ be a vector space over $F$, $\bold{t}$ be a (multi) parameter.
A field $a(z)$ (on $V$) is a series of the form 
\begin{displaymath}
a(z)=\sum_{n\in \mathbf{Z}}a_{(n)}z^{-n-1}, \ \ \ a_{(n)}\in
\End(V),\ a_{(n)}v=0 \ \text{for any}\ v\in V, \ n\gg 0.
\end{displaymath}
A state-field correspondence is a linear  map from $V$ to the space of fields
 that associates to any $a\in V$ a field $Y(a,z)=a(z)$.
We will use also the following notation introduced in \cite{EK}
\begin{align*}
& Y: V\ten V \to V((z)),\\ & a\ten b \mapsto Y(z)(a\ten b).
\end{align*}
Denote by $\mathbf{K}$ the algebra \mbox{$k[z^{\pm 1}, w^{\pm 1}, (z-w)^{\pm 1}][[\bold{t}]]$.}
\begin{defn}
An  $H_D$-quantum vertex algebra  over $k[[\bold{t}]]$ consists of the following data:
\begin{itemize}
\item the space of states---a topologically free module $V$ over $k[[\bold{t}]]$.
\item the vacuum vector---a vector $|0\rangle  \in V$.
\item the space of fields and state-field correspondence.
\item a distinguished operator $D: V\to V$
\item a braiding map---a linear map \mbox{$R(z,w):V\ten V \to V\ten V \ten \mathbf{K}$}, which satisfies the conditions:
\begin{itemize}
\item   Yang-Baxter equation: 
\begin{equation*} 
R_{12}(z_1, z_2)R_{13}(z_1, z_3)R_{23}(z_2, z_3)=R_{23}(z_2, z_3)R_{13}(z_1, z_3)R_{12}(z_1, z_2)
\end{equation*}
\item shift conditions
\begin{align*}
& [D\ten 1, R(z,w)]=-\partial _zR(z,w),\\ & [1\ten D, R(z,w)]=-\partial _wR(z,w).
\end{align*}
\item unitarity condition
\begin{equation*}
\tau R_{w,z}\tau =R_{z,w}^{-1}, \ \text{where} \  \tau  \ \text{is the flip,} \ \tau (a\ten b)=b\ten a.
\end{equation*}
Here  we denote $R_{z,w}=R(z,w)$.
\end{itemize}
\end{itemize}
These data should satisfy the following  two sets of axioms 
\begin{itemize}
\item The axioms related to the action of the Hopf algebra
  $H_D=\mathbf{C} [D]$: 
\begin{itemize}
\item translation covariance:\ \  $Y(Da,z)=\partial _zY(a,z)$.
\item vacuum axiom: \ \ $Y(|0\rangle ,z)=Id_V$.
\item creation axiom: \ \ $Y(a,z)|0\rangle \arrowvert _{z=0} =a$.
\end{itemize}
\end{itemize}
\begin{itemize}
\item Braided locality axiom: for all $a,b \in V$  there exist $N$ such that 
\begin{equation*}
(z-w)^N Y(a,z)Y(b,w)=(w-z)^N \tilde{R}(Y(b,w),Y(a,z)), 
\end{equation*}
where we denote  
\begin{equation}
\label{eq:rtilde}
 \tilde{R}(Y(b,w),Y(a,z))c=Y(w)(1\ten Y(z))(R(w,z)(b\ten a)\ten c), 
\end{equation}
for any $a,b,c \in V$.
\end{itemize}
\end{defn}
This definition is a generalization of the definition from  \cite{EK}
and therefore we will sometimes use the name  vertex algebras of EK type for
$H_D$-quantum  vertex algebras. In order to accommodate the examples of the
symmetric polynomials the following two important changes were made.
First, in \cite{EK} the braiding $R(z,w)$ is in fact of the form
$\hat{R}(z-w)$, where $\hat{R}$ is a function of a single
variable. This is no longer possible for quantum vertex algebras
describing the symmetric polynomials, where the braiding map is truly
a function of two variables (see the appendix in \cite{A2}). Second,
in \cite{EK} (as in the case of classical vertex algebras) the vertex
operators satisfy 
\begin{equation}
\label{eq:wrongtrcov}
\frac{d}{dz}Y(z)=DY(z)-Y(z)(1\ten D).
\end{equation}
The above equation is not valid
for Hall-Littlewood vertex operators for instance, but the translation
covariance as required in the axiom holds (see the
appendix in \cite{A2}, as well as \cite{AB}). This last consideration in fact implies that our $H_D$-quantum  vertex
algebras are no longer field algebras (as defined for instance in \cite{Kac}).

The classical vertex (super)algebras are a particular case of
$H_D$-quantum vertex algebras, with braiding map
$R(z,w)=\pm Id_V\ten Id_V$, and for them \eqref{eq:wrongtrcov} is
equivalent to the translation covariance axiom. 
In order to prove that the quantum vertex operators associated with the symmetric polynomials belong to the type of structure described above we have to follow a somewhat roundabout way.

\section{Borcherds bicharacter construction of vertex algebras}

This section recalls Borcherds bicharacter construction from
\cite{Borc}. This construction  was used there to construct examples
of  $(A,H,S)$ (quantum) vertex algebra, but we will not discuss here
$(A,H,S)$ vertex algebras, as it is outside of the scope of this
paper. We will adapt the bicharacter construction as a tool for
recovering the  braiding map between the quantum vertex operators. The idea is
to define a braided vertex algebra using a bicharacter (it can be made
into $(A,H,S)$ quantum vertex algebra, see Theorem 4.2, \cite{Borc}),
as such algebras come with a formula for the braiding map, then to prove that the vertex operators  we are dealing with can be  in fact identified with those defined via a bicharacter. We will provide  in \cite{AB} some of the details which are omitted here.

First let us recall some definitions.
\subsection{Bicharacters}
\label{bicharacters}

Let $M$ be a commutative and cocommutative Hopf algebra. 
Denote the coproduct and the counit by $\del$  and $\eta$, the antipode by $S$

If $a$ is an element of a coalgebra  we will use Sweedler's notation and write $\del (a)=\Sigma \ a^{\prime}\ten a^{\second}$. We often will omit the summation sign.
\begin{defn}(Bicharacter) A bicharacter on $M$ is
a linear map $r$ from $M\ten M$ to $T$ , where $T$ is a commutative
$k$ algebra, such that for any $a,b,c\in M$
\begin{align*}
r(1\ten a)&= \eta (a) = r(a\ten 1),\\
r(ab\ten c)&= \Sigma r(a\ten c^{\prime})r(b\ten c^{\second}),\\
r(a\ten bc)&= \Sigma r(a^{\prime}\ten b)(r(a^{\second}\ten c). 
\end{align*}
\end{defn}
We will be interested mostly in the cases where the target space for the
bicharacter is $\mathbf{K}$. 

Note that in fact if $T=\mathbf{K}$ the bicharacter $r(a\ten b)$ is a function of $z$ and $w$, and should be written as $r(a\ten b)(z,w)$ but we will omit the $(z,w)$  unless strictly necessary, in which case we will sometimes write it as $r(a_z\ten b_w)$. 

\begin{lem/defn}
Let $r$ and $s$ be bicharacters, where $M$ is a commutative cocommutative
Hopf algebra. We can define a convolution product $r*s$ by 
\begin{equation*}
r*s (a\ten b)= r(a^{\prime}\ten b^{\prime})s(a^{\second}\ten b^{\second}). 
\end{equation*}
The identity bicharacter is given by \mbox{$\epsilon(a\ten b)= \eta (a) \eta
(b) $.}\\
The inverse bicharacter is defined by \mbox{$r^{-1}(a\ten b)=r(S(a)\ten b)$.}\\
The transpose bicharacter is defined to be \mbox{$r^{\tau }(a_z\ten b_w)=r(b_w\ten a_z)$.} (Here we assume that  the target space $T=\mathbf{K}$.)\\
The bicharacters on $M$ form a commutative group.
\end{lem/defn} Note: In the above definition it is essential that  $M$ is
a cocommutative Hopf algebra.

\subsection{Free Leibnitz modules and extension of bicharacters}
\label{sec:free module}
  \ \ 

Now we proceed to describe the type of vector spaces  we will be using as underlying our braided vertex algebras. 
Denote by $H_D$  the Hopf algebra  $\mathbf{C}[D]$, generated by a primitive element $D$.
As usual we write $D^{(n)}=\frac{D^n}{n!}$. We have then
\begin{align}
& \del D^{(n)}=\sum _{k+l=n} D^{(k)}\ten D^{(l)},\ \  S(D^{(n)})=(-1)^nD^{(n)}.
\end{align}
For any $M$ a commutative, cocommutative Hopf algebra we want to
construct a (universal) $H_D$-module algebra $H_D(M)$ containing
$M$.\\ Let \mbox{$H_D(M)\equiv \Sym (H_D\ten M)$.} Then $H_D(M)$ is
naturally a module-algebra over $H_D$. We extend the
comultiplication, the counit and the antipode from $M$ to $H_D(M)$ as follows.
We have for any $a\in H(M), \ h\in H_D$ 
\begin{align*}
& \del h.a=\sum h^{\prime}.a^{\prime}\ten h^{\second}.a^{\second},\\
& \eta(h.a)=\eta (h)(\eta (a)),\\
& S(h.a)=S(h).S(a).
\end{align*}
It is easy to check that the comultiplication, the counit and the
antipode defined as above will turn $H_D(M)$ into commutative
cocommutative Hopf algebra.
All of the vector spaces underlying our vertex algebras in this paper
are going to be of the type $H_D(M)$ for some $M$. We call such
module-algebras free $H_D$-Leibnitz modules.

If we have a bicharacter on $M$ with target $T$ which is an $H_D$-bimodule we can extend it to $H_D(M)$ by requiring
\begin{equation*}
r(Da\ten b)(z,w)=\partial _z r(a\ten b)(z,w), \ \ 
r(a\ten Db)(z,w)=\partial _w r(a\ten b)(z,w).
\end{equation*}
It is easy to check that this extension satisfies the axioms for a bicharacter on $H_D(M)$.

\subsection{Fields defined via bicharacter}
\label{sec:fields}
 \ \ \ 

Let $V=H_D(M)$ for some commutative cocommutative algebra $M$ with a bicharacter $r$  with target $\mathbf{K}$.  
\begin{defn}
For any $a,b\in V$ define $Y(a,z)b\in V((z))$ by
\begin{equation}
\label{eq:fields}
Y(a,z)b=\sum (e^{zD}(a^{\prime}))b^{\prime}r(a^{\second}\ten b^{\second})(z,0)
\end{equation}
\end{defn}
This definition gives us a  state-field correspondence, even though  the fields are not given by generating series as is usual.
Since  $M$ and hence  $V$ is an algebra with unit element $1$,
define the vacuum vector for the quantum vertex algebra structure we are constructing to be the unit
element in the algebra $V$. We have the following lemma:
\begin{lem}
The fields $Y(a,z)b$ for $a,b\in V$ defined by \eqref{eq:fields}
satisfy the vacuum, the creation and the translation covariance axioms
for an $H_D$-quantum  vertex algebra.
\end{lem}
\begin{proof}
For all $a\in V$
\begin{align*}
& Y(a,z)1=\sum e^{zD}(a^{\prime})r(a^{\second}\ten 1)(z,0)=\\
& =\sum e^{zD}(a^{\prime})\eta (a^{\second})=\sum a^{\prime}\eta
(a^{\second})+z.O(z)=a+z.O(z), 
\end{align*}
which proves  the creation axiom.
Similarly
\begin{align*}
& Y(Da,z)b=\sum e^{zD}((Da)^{\prime})b^{\prime}r((Da)^{\second}\ten b^{\second})(z,0)=\\
& = \sum e^{zD}(Da^{\prime})b^{\prime}r(a^{\second}\ten b^{\second})(z,0)+\sum
e^{zD}(a)^{\prime}b^{\prime}r(Da^{\second}\ten b^{\second})(z,0)=\\ 
&=\partial _z (\sum e^{zD}a^{\prime}b^{\prime}r(a^{\second}\ten
b^{\second})(z,0))=\partial _zY(a,z)b. 
\end{align*}
This proves the  translation covariance, the vacuum axiom is proved
similarly.
\end{proof}
\begin{defn}
Let $a,b\in V$. For given bicharacter we define 
\begin{equation*}
R(z,w): V\ten V \to V\ten V \ten  \mathbf{K} 
\end{equation*}
by
\begin{equation*}
R(a\ten b)=\Sigma a^{\prime}\ten b^{\prime}r^t*r^{-1}(a^{\second}\ten b^{\second})
\end{equation*}
\end{defn}

The map $R$ is going to be our braiding map. It is easy to prove that
it satisfies the Yang-Baxter equation, the shift  and the unitarity
conditions (more details can be found in \cite{AB}).
We have the following 
\begin{thm}
\label{thm:axioms}
The  data defined as above (the space of states $V$, the vacuum
vector $1$, the distinguished operator $D$, the state-field
correspondence as defined by \eqref{eq:fields} above and  the braiding
map $R$) satisfies the axioms of $H_D$-quantum vertex algebra.
\end{thm}

\section{Bicharacter construction for lattice (quantum) vertex algebras}
\label{rank-one lattice}
In this section we connect the vertex operators as defined in the previous section to the vertex operator $\Psi$ defined in the first section.
We consider a specific case of the very general bicharacter construction described in the previous section.

Let $M= \mathbf{C}[\mathbf{Z}\alpha]$,  the group algebra of the
rank-one lattice generated by $e^{m\alpha}$ (such that
\mbox{$e^{m\alpha}e^{n\alpha }=e^{(m+n)\alpha}$).} Then $V=H_D (M)$ can be written as \\ \mbox{$\mathbf{C}[h^{(n)}]\ten
\mathbf{C}[\mathbf{Z}\alpha]$, $n\ge 1$,}  where  the new variables
$h^{(n)}$ are defined by 
\begin{equation}
h=D\alpha=(D.e^{\alpha })e^{-\alpha },\ \ h^{(n)}=D^{(n)}h.
\end{equation}
The bialgebra structure on $V$ is determined by $e^{\alpha}$ being grouplike. For instance we have 
\begin{equation*}
\del(h)=h\ten 1 + 1\ten h,
\end{equation*}
i.e.,  $h$ is a primitive element.
The antipode is determined by $S(e^{\alpha})=e^{-\alpha}$. It follows then that \begin{equation*}
S(h)=-h.
\end{equation*}
To define a bicharacter on $V$ we only need to define it on $M$. 
Moreover,
since \mbox{$r(e^{\alpha}\ten e^{-\alpha})r(e^{\alpha}\ten e^{\alpha})=1$}
(follows directly from the properties of a bicharacter), 
we only need to define \mbox{$\sigma (z,w)=r(e^{\alpha}\ten e^{\alpha})(z,w)$.}
Simple calculation shows that 
\begin{align*}
& r(e^{\alpha }\ten h)(z,w)=\frac{\partial _w \sigma (z,w)}{\sigma (z,w)}=\partial _w \ln\sigma (z,w),\\
& r(h\ten h)(z,w)=\partial _z \partial _w \ln\sigma (z,w).
\end{align*}

For any function $f(x; \bold{t})$ depending on a (multi) parameter denote by
$i_{\bold{t}}f(x; \bold{t})$ the  power series expansion of $f(x; \bold{t})$ in
terms of this parameter (i.e. consider the parameter(s) as very small). 
We will suppress the dependence on the parameter
in the notation and will write just $f(x)$ instead of $f(x, \bold{t})$.

Let 
\begin{align}
\label{eq:f}
&f(x)=\exp (-\sum _{n\ge 1} \frac{v^{-1}_n}{n} x^n ),\\
\label{eq:sigma}
&\sigma (z,w)=z.i_{\bold{t}}f(\frac{w}{z}).
\end{align}
The choice of the (multi)parameter $\bold{t}$ in the above expansion will depend
on the multiplicative family $(v_{\lambda })$.
For instance, in the case of Hall-Littlewood polynomials we have 
\begin{equation}
\label{eq:param}
\bold{t}=t, \ \ v_n^{-1}=1-t^n, \  
f(x)=\frac{1-x}{1-tx}, \ \ \ i_{\bold{t}}f(x)=(1-x)\sum_{n=0}^{\infty}x^nt^n.
\end{equation}

We have the following two lemmas which provide the connection between
the bicharacter definition of fields and the expressions as 
generating series as for instance in  Lemma \eqref{eq:psi}.
\begin{lem}
Let $\sigma (z,w)$ be as above. The field \mbox{$h(z):=Y(D\alpha ,z)$}
defined by \eqref{eq:fields} is a Heisenberg field, i.e.,  if we write
$h(z)$ as \mbox{$h(z)=\sum _{n\in \mathbf{Z}}h_n z^{-n-1}$,} then the
operators \mbox{$h_n\in End(V),\ n\in \mathbf{Z}$} define a deformed Heisenberg algebra (as in Lemma \ref{lem:heis}).
\end{lem}
\begin{lem}
\label{lem:expon}
Let $\sigma (z,w)$ be as above.  The field $Y(e^{\alpha },z)$ as defined in \eqref{eq:fields} can be written as 
\begin{equation*}
Y(e^{\alpha },z)=\exp (\sum _{n\ge 1}\frac{h_{-n}}{n} z^n)\exp (-\sum _{n\ge 1}\frac{h_{n}}{n} z^{-n})e^{\alpha}z^{\partial _{\alpha }},
\end{equation*}
where the operators $h_n$ are the generators of the Heisenberg algebra.
\end{lem}
For proofs, see \cite{A2} or \cite{AB}.
Hence we have the following theorem:
\begin{thm}
\label{thm:main}
Let $\sigma (z,w)$ be defined as in \eqref{eq:sigma}. The field $Y(e^{\alpha
},z)$ as defined in \eqref{eq:fields} equals the field $\Psi(z)$ as
defined in \eqref{eq:psi}.Moreover, the field $\Psi(z)$ 
belongs to an $H_D$-quantum vertex algebra structure on \mbox{$V=\Sym(D^{(n)}\alpha)\ten \mathbf{C}[\mathbf{Z}\alpha]$,} with braiding map given by
\begin{equation}
\label{eq:braiding}
R(a\ten b)=\Sigma a^{\prime}\ten b^{\prime}r^t*r^{-1}(a^{\second}\ten b^{\second}),
\end{equation}
for any $a,b\in V$.
\end{thm}
For instance, in the simplest case of $a=e^{\alpha },b=e^{\alpha }$,
the braiding map looks as follows:
\begin{equation}
\label{eq:simple}
\tilde{R}(Y(e^{\alpha },z),Y(e^{\alpha },w))=i_{\bold{t}}\Big(
\frac{w.f(\frac{z}{w})}{z.f(\frac{w}{z})}\Big) Y(e^{\alpha },z)Y(e^{\alpha },w)
\end{equation}
where $\tilde{R}$ is defined as in \eqref{eq:rtilde}.

The function $f(x)$ as defined in \eqref{eq:f} is precisely the
one considered by Jing in \cite{J3}, where the above equation
\eqref{eq:simple} is also derived. The use of $i_{\bold{t}}f(x)$ instead of $f(x)$ is dictated by
the choice of $\mathbf{K}$ as the target space for the bicharacters, and
that in turn is chosen so that we can view our quantum vertex algebras as a
kind of ``deformation quantization'' of the classical ones. It is
notable that although the braiding map $R(z,w)$ is an expansion in power
series in the small parameter, it is also a double-sided series in $z$ and
$w$. Consider for instance  in the case of Hall-Littlewood
polynomials (see \eqref{eq:param}) the braiding $R(h\ten h)(z,w)$ (where $h=D\alpha $). We have from \eqref{eq:braiding}
\begin{align*}
R(h\ten h)(z,w)&=h\ten h +1\ten 1i_t\Big(-
\frac{t}{(z-tw)^2}+\frac{t}{(w-zt)^2}\Big)=\\
&=h\ten h +1\ten 1\Big(-i_{z,w}\frac{t}{(z-tw)^2}+i_{w,z}\frac{t}{(w-zt)^2}\Big),
\end{align*}
where for a meromorphic function $f(z,w)$ we denote by $i_{z,w}f(z,w)$
the expansion of $f(z,w)$ in the region $z\gg w$, and similarly for
$i_{w,z}f(z,w)$.  More details are provided in \cite{AB}.

For the case of Hall-Littlewood polynomials  the vertex operator $\Psi _t$ 
\begin{equation*}
\Psi _t(z)=\exp (\sum
_{n\ge 1} \frac{1-t^n}{n} p_n z^n)\exp (-\sum _{n\ge
  1}\frac{1-t^n}{n} p^{\perp}_{n} z^{-n})e^{\alpha}z^{\partial
  _{\alpha}}
\end{equation*}
belongs to an $H_D$-quantum vertex algebra over
$k[[t]]$. If $t=0$ we get back  the classical vertex superalgebra
describing the Schur polynomials.

In the case corresponding to Macdonald's polynomials we have 
\begin{equation*}
\bold{t}=(t,q), \ \ v_n^{-1}=\frac{1-t^n}{1-q^n}, \  \ f(x)=\frac{(x;q)_{\infty}}{(tx;q)_{\infty}},
\end{equation*}
where $(x;q)_{\infty}=\prod_{n\ge 0}(1-xq^n)$. The vertex operator
\begin{align*}
\Psi _{q,t}(z)=\exp (\sum
_{n\ge 1} \frac{1-t^n}{n(1-q^n)} p_n z^n)\exp (-\sum _{n\ge
  1}\frac{1-t^n}{n(1-q^n)} p^{\perp}_{n}z^{-n})e^{\alpha}z^{\partial
  _{\alpha}}
\end{align*}
belongs to $H_D$-quantum vertex algebra structure over $k[[t,q]]$. If
$q=0$ we get back the operator $\Psi _t$.

We want to make a note that all the bosonic vertex operators (i.e., vertex operators of exponential form  similar to $\Psi(z)$) can be written in  terms of bicharacters. For instance the bosonic vertex representation of the quantum affine algebras (as considered for instance by I. Frenkel and Jing in \cite{FJ}) can be written using a bicharacter. That is a help toward defining the vertex algebra structure, as we have a formula for the braiding map, but is only part of the picture, because as it turns out these vertex  representations  cannot be accommodated by the axioms of $H_D$-quantum  vertex algebras. 
The $H_D$-quantum  vertex
algebras are part of bigger structures---FR type of deformed chiral
algebras, that involve bigger Hopf algebras, containing $H_D$. This
will be discussed in \cite{A2}, where we also
extend the definition of ``fields'' so that we can incorporate  more
general Hopf algebras. Thus the $H_D$-quantum vertex algebras are, if not the complete answer, at least a necessary ingredient in the description of quantum  vertex algebras.

\end{document}